\documentclass[12pt, reqno]{amsart}
\usepackage{ amsmath,amsthm, amscd, amsfonts, amssymb, graphicx, color}
\usepackage[bookmarksnumbered, colorlinks, plainpages]{hyperref}
\newtheorem{thm}{Theorem}[section]
\newtheorem{cor}[thm]{Corollary}
\newtheorem{lem}[thm]{Lemma}
\newtheorem{prop}[thm]{Proposition}

\newtheorem{exam}[thm]{Example}
\numberwithin{equation}{section}

\begin{document}

\title{Hirano inverse of anti-triangular matrix over Banach Algebras}

\author[Gou Haibo]{Gou Haibo}
\address{Gou Haibo\\ School of Mathematics \\ Hangzhou Normal University \\ China}
\email{<haibo\_gou@163.com>}

\author[Chen Huanyin]{Chen Huanyin}
\address{Chen Huanyin \\ School of Mathematics \\ Hangzhou Normal University \\ China}
\email{<huanyinchen@aliyun.com>}
\thanks{}


\begin{abstract}
In this paper we investigate Hirano invertibility of anti-triangular matrix over a Banach algebra. Let $a\in {\mathcal A}^H, b\in {\mathcal A}^{sD}.$ If $b^Da=0, bab^{\pi}=0,$ we prove that $\begin{pmatrix}
a&1\\
b&0
\end{pmatrix}\in M_2(\mathcal A)^H.$ Moreover, we considered Hirano invertibility of anti-triangular matrices under commutative-like conditions. These provide new kind of operator matrices with tripotent and nilpotent decompositions.
\end{abstract}

\maketitle

\section{Introduction}

Hirano invertibility was introduced by Chen and Abdolyousefi (see [10]). An element $a\in \mathcal A$ has Hirano inverse if there exists $x\in \mathcal A$ such that $$ax=xa, x=xax, a^2-ax\in \mathcal A$$ is nilpotent. Such $x,$ if it exists, is unique and will be denoted by $a^H.$

Recall that an element $a\in \mathcal A$ has strongly Drazin inverse if there exists $x\in \mathcal A$ such that $$ax=xa, x=xax, a-ax\in \mathcal A$$ is nilpotent. Such $x,$ if it exists, is unique and will be denoted by $a^{sD}.$
Evidently, $a\in \mathcal{A}^{sD}$ if and only $a-a^2\in \mathcal{A}$ is nilpotent if and only if $a$ is the sum of an idempotent and a nilpotent  that commute (see [8, Lemma 2.2]).  In this paper, we focus on the Hirano invertibility of anti-triangular matrices over a Banach Algebra.

These provide new kind of operator matrices which are the sum of an tripotent and a nilpotent.

In Section 2, we present several elementary results of Hirano inverse, which will be repeatedly used in the sequel.

In Section 3, we study the Hirano invertibility of anti-triangular operator matrices with a form $\begin{pmatrix}
a&1\\
b&0
\end{pmatrix}$ when $a\in {\mathcal A}^H, b\in {\mathcal A}^{sD}, b^Da=0$ and $bab^{\pi}=0.$

In Section 4, we consider some communicative-like conditions on anti-triangular matrix to study its Hirano invertibility. Let $a\in {\mathcal A}^H, b\in {\mathcal A}^{sD}.$ If $b^Da=0, bab^{\pi}=abb^{\pi},$ we prove that $\begin{pmatrix}
a&1\\
b&0
\end{pmatrix}\in M_2(\mathcal A)^H.$

In Section 5, several examples are given to illustrate our results.

Throughout the paper, $\mathcal A$ denotes a complex Banach algebra with an identity. $\mathcal A^{H}$ represents the set of all Hirano invertible elements in $\mathcal A.$ $\mathcal A^{sD}$ is the set of all elements which are strongly Drazin invertible. $N(\mathcal A)$ is the set of nilpotent elements in $\mathcal A$ and $U(\mathcal A)$ stands for the set of all invertible elements in $\mathcal A.$

\section{Key lemmas}

The aim of this section is to present some elementary results of Hirano inverse which will be used  in the sequel.

\begin{prop}
Let $\mathcal A$ be a Banach algebra. Then the following are equivalent:
\begin{enumerate}
\item $a\in \mathcal A^{H}.$
\item $a-a^3\in \mathcal A$ is nilpotent.
\item $a$ is the sum of a tripotent and a nilpotent that commute.
\end{enumerate}
\end{prop}

\begin{proof} 
See [7, Theorem 3.1] and [11, Corollary 2.8].

 We have (1)$\Leftrightarrow$(2)$\Leftrightarrow$(3) obtained.
\end{proof}

\begin{lem}
If $a,b\in \mathcal A^{H},$ then $\begin{pmatrix}
a&c\\
0&b
\end{pmatrix} \in M_2(\mathcal A)^{H}.$
\end{lem}

\begin{proof} 
Since $a,b\in \mathcal A^H,$ we have $(a-a^3)^m=0$ and $(b-b^3)^n=0$ for some $m,n\in \mathbb N.$
According to Proposition 2.1, we only need to prove that 
$$
\begin{pmatrix}
a&c\\
0&b
\end{pmatrix}-\begin{pmatrix}
a&c\\
0&b
\end{pmatrix}^3\in N(\mathcal A).
$$ 
$$\begin{pmatrix}
a&c\\
0&b
\end{pmatrix}-\begin{pmatrix}
a&c\\
0&b
\end{pmatrix}^3=\begin{pmatrix}
a-a^3&c-(a^2c+abc+cb^2)\\
0&b-b^3
\end{pmatrix}.$$ Let $A=a-a^3, B=c-(a^2c+abc+cb^2), D=b-b^3, X=\begin{pmatrix}
A&B\\
0&D
\end{pmatrix}.$ We directly compute that
$$
\begin{array}{r l l}
X^{n+m}&=&{\begin{pmatrix}
A^{(n+m)}&\sum\limits_{i=0}^{(n+m)-1}A^{(n+m)-1-i}BD^i\\
0&D^{(n+m)}
\end{pmatrix}}\\
&=&{\begin{pmatrix}
0&\sum\limits_{i=0}^{(n+m)-1}A^{(n+m)-1-i}BD^i\\
0&0
\end{pmatrix}}.
\end{array}
$$
Then we have $X^{2(n+m)}=0,$ which infers that
$$
\begin{pmatrix}
a&c\\
0&b
\end{pmatrix}-\begin{pmatrix}
a&c\\
0&b
\end{pmatrix}^3\in N(\mathcal A).
$$
Therefore, by Proposition 2.1, $\begin{pmatrix}
a&c\\
0&b
\end{pmatrix}\in M_2(\mathcal A)^H.$
\end{proof}

\begin{lem}
Let $A\in M_{m\times n}(\mathcal A), B\in M_{n\times m}(\mathcal A).$ If $AB\in M_m(\mathcal A)^{H},$ then $BA\in M_n(\mathcal A)^{H}.$
\end{lem}

\begin{proof} 
We may assume that $m\geq n.$ Then $AB=\begin{pmatrix}
A&0\\
\end{pmatrix}\begin{pmatrix}
B\\0
\end{pmatrix}.$ In light of [10],  we have $\begin{pmatrix}
B\\0
\end{pmatrix}\begin{pmatrix}
A&0\\
\end{pmatrix}$ has Hirano inverse. That is $\begin{pmatrix}
BA&0\\
0&0
\end{pmatrix}$ has Hirano inverse. In view of Proposition 2.1, we have $\bigg{[}\begin{pmatrix}
BA&0\\
0&0
\end{pmatrix}-{\begin{pmatrix}
BA&0\\
0&0
\end{pmatrix}}^3\bigg{]}^m=0.$ Therefore, $[BA-(BA)^3]^m=0$ Accordingly, $BA\in M_n(\mathcal A)^{H},$ as asserted.
\end{proof}

\begin{lem}
Let $a, b\in \mathcal A^{H}$. If $ab=0,$ then $a+b\in \mathcal A^{H}.$
\end{lem}

\begin{proof}
Let $A=\begin{pmatrix}
1\\a
\end{pmatrix}, B=\begin{pmatrix}
b&1
\end{pmatrix}$, then $a+b=BA, AB=\begin{pmatrix}
b&1\\
0&a
\end{pmatrix}.$ According to Lemma 2.3, $(BA)^{H}=B{[(AB)^{H}]^2}A.$ By using of Lemma 2.2 and Lemma 2.3, $AB\in M_2(\mathcal A)^{H}.$ Therefore $BA=a+b\in \mathcal A^{H}.$
\end{proof}

\begin{thm} 
Let $a, b\in {\mathcal A}^H.$ If $aba=0$ and $ab^2=0,$ then $a+b\in {\mathcal A}^H.$
\end{thm}

\begin{proof} 
Let $a+b=\begin{pmatrix}
1&b
\end{pmatrix}\begin{pmatrix}
a\\
1
\end{pmatrix}.$ Then $a+b\in {\mathcal A}^H$ $\Leftrightarrow$ 
$\begin{pmatrix}
a\\
1
\end{pmatrix}\begin{pmatrix}
1&b
\end{pmatrix}=\begin{pmatrix}
a&ab\\
1&b
\end{pmatrix}$
$\in {\mathcal A}^H$ according to Lemma 2.3. We deduce that

$$\begin{pmatrix}
a&ab\\
1&b
\end{pmatrix}=\begin{pmatrix}
a&0\\
0&0
\end{pmatrix}+\begin{pmatrix}
0&0\\
0&b
\end{pmatrix}+\begin{pmatrix}
0&ab\\
1&0
\end{pmatrix}.$$
Since $a, b\in {\mathcal A}^H,$ then $\begin{pmatrix}
a&0\\
0&0
\end{pmatrix}$ and $\begin{pmatrix}
0&0\\
0&b
\end{pmatrix}\in M_2(\mathcal A)^H.$
$$
\begin{pmatrix}
0&ab\\
1&0
\end{pmatrix}-{\begin{pmatrix}
0&ab\\
1&0
\end{pmatrix}}^3=\begin{pmatrix}
0&ab\\
1-ab&0
\end{pmatrix},
$$
$$
{\begin{pmatrix}
0&ab\\
1-ab&0
\end{pmatrix}}^4=0.
$$
According to Proposition 2.1, $\begin{pmatrix}
0&ab\\
1&0
\end{pmatrix}\in M_2(\mathcal A)^H$ $\Rightarrow$ $\begin{pmatrix}
a&ab\\
1&b
\end{pmatrix}\in M_2(\mathcal A)^H$ $\Rightarrow$ $a+b\in {\mathcal A}^H$ obtained.
\end{proof}

\section{Anti-triangular operator matrices}

The purpose of this section is to investigate the Hirano invertibility of anti-triangular operator matrices. We now derive them.

\begin{thm}
Let $x=\begin{pmatrix}
a&1\\
b&0
\end{pmatrix},$
$a\in {\mathcal A}^H, b\in \mathcal A^{sD} .$ If $b^Da=0, bab^{\pi}=0$, then $x\in M_2(\mathcal A)^{H}.$
\end{thm}

\begin{proof}
Let $e=\begin{pmatrix}
bb^{D}&0\\
0&1
\end{pmatrix},$ then $e^2=e.$ Set $bb^{D}=b^e.$ Then we have the Peirce decomposition of $x$ related to the idempotent $e$.
$$x=\begin{pmatrix}
exe&ex(1-e)\\
(1-e)xe&(1-e)x(1-e)
\end{pmatrix}_e=\begin{pmatrix}
\alpha&\beta\\
\gamma&\delta
\end{pmatrix}_e=\alpha+\beta+\gamma+\delta,$$
We compute that
 $$
\begin{array}{r l l}
\alpha&=&\begin{pmatrix}
bb^D&0\\
0&1
\end{pmatrix}\begin{pmatrix}
a&1\\
b&0
\end{pmatrix}\begin{pmatrix}
bb^D&0\\
0&1
\end{pmatrix}\\
&=&\begin{pmatrix}
bb^Dabb^D&bb^D\\
bbb^D&0
\end{pmatrix}=\begin{pmatrix}
bb^Da&bb^D\\
b^2b^D&0
\end{pmatrix}\\
&=&
\begin{pmatrix}
0&b^e\\
bb^e&0
\end{pmatrix}_e,
\end{array}
$$
$$
\beta=\begin{pmatrix}
0&0\\
bb^{\pi}&0
\end{pmatrix}_e,
\gamma=\begin{pmatrix}
b^{\pi}ab^e&b^{\pi}\\
0&0
\end{pmatrix}_e,
\delta=\begin{pmatrix}
ab^{\pi}&0\\0&0
\end{pmatrix}_e, $$
where $\alpha, \beta, \gamma, \delta \in \mathcal A.$ We compute that
$${\alpha}^3=\begin{pmatrix}
0&b^2b^{D}\\
b^3b^{D}&0
\end{pmatrix},$$
$$\alpha-{\alpha}^3=\begin{pmatrix}
0&(b-b^2)b^{D}\\
(b-b^2)bb^{D}&0
\end{pmatrix},$$ By hypothesis, $b$ has strongly Drazin inverse, and so $b-b^2\in \mathcal{A}$ is nilpotent.
Thus $\alpha$ is Hirano invertible by Proposition 2.1.
Note that $(\beta +\gamma)^2=\begin{pmatrix}
bb^{\pi}&0\\
bb^{\pi}a&bb^{\pi}
\end{pmatrix}$ is nilpotent, and then $\beta +\gamma$ is nilpotent too. Thus, $(\beta +\gamma)^{H}=0.$
Obviously, we have $$ab^{\pi}-(ab^{\pi})^3=ab^{\pi}-ab^{\pi}ab^{\pi}ab^{\pi}=(a-a^3)b^{\pi},$$ and then
$$[ab^{\pi}-(ab^{\pi})^3]^m=[(a-a^3)b^{\pi}]^m=\cdots=(a-a^3)^mb^{\pi}$$ for any $m\in {\Bbb N}$. This implies that $ab^{\pi}\in {\mathcal A}^H.$
Thus $\delta$ is Hirano invertible. For $\beta\delta=0, \gamma\delta=0,$ so $(\beta +\gamma)\delta=0.$ According to Lemma 2.4, $ \beta+\gamma+\delta \in \mathcal A^{H}.$\\
$$\beta+\gamma+\delta=x-\alpha=
\begin{pmatrix}
a&1\\
b&0
\end{pmatrix}-\alpha=
\begin{pmatrix}
b^{\pi}a&b^{\pi}\\
bb^{\pi}&0
\end{pmatrix},$$ and we have $\alpha(\beta+\gamma+\delta)=0.$ According to Lemma 2.4, $ \alpha+\beta+\gamma+\delta \in \mathcal A^{H}.$ Therefore, $x\in M_2(\mathcal A)^{H}.$
\end{proof}

\begin{cor}
Let $x=\begin{pmatrix}
a&b\\
c&d
\end{pmatrix},$
$a,d\in {\mathcal A}^H, bc\in \mathcal A^{sD} .$ If $(bc)^Da=0, bca(bc)^{\pi}=0, bdc=0$ and $bd^2=0$, then $x\in M_2(\mathcal A)^{H}.$
\end{cor}

\begin{proof}
Write $x=y+z,$ where $y=\begin{pmatrix}
0&0\\
0&d
\end{pmatrix}$ and $z=\begin{pmatrix}
a&b\\
c&0
\end{pmatrix}.$
Obviously, $y$ is Hirano invertible and $y^{H}=
\begin{pmatrix}
0&0\\
0&d^{H}
\end{pmatrix}.$ To prove that $z$ is Hirano invertible, we write $z=pq,$ where $p=
\begin{pmatrix}
a&1\\
c&0
\end{pmatrix}$ and $q=
\begin{pmatrix}
1&0\\
0&b
\end{pmatrix}.$ Applying Theorem 3.1, we deduce that $qp=
\begin{pmatrix}
a&1\\
bc&0
\end{pmatrix}$ is Hirano invertible. According to Lemma 2.3, $z=pq$ is Hirano invertible and $z^{H}=p[(qp)^{H}]^2q.$
Since $bdc=0$ and $bd^2=0,$ we see that $zyz=0, zy^2=0.$ By using of Theorem 2.5, we have $x=y+z\in M_2(\mathcal A)^{H}.$
\end{proof}

\begin{cor}
Let $x=\begin{pmatrix}
a&b\\
c&d
\end{pmatrix},$
$a,d\in {\mathcal A}^H, bc\in \mathcal A^{sD}.$ If $(bc)^Da=0, a{(bc)}^{\pi}=0$ and $bd=0,$ then $x\in M_2(\mathcal A)^{H}.$
\end{cor}

\begin{proof}
This is an immediate application of Corollary 3.2.
\end{proof}

\begin{lem}
Let $C\in \mathbb C^{n\times n}.$ Then $C\in (\mathbb C^{n\times n})^H$ if and only if the eigenvalues of $C$ are $-1$, $0$, $1$.
\end{lem}

\begin{proof}
If $C\in {\Bbb C}^{n\times n}$ is Hirano invertible, then there exists $P,$ such that $$P^{-1}CP=\begin{pmatrix}
J_1&0&\cdots&0\\
0&J_2&\cdots&0\\
\vdots&\vdots&\ddots&\vdots\\
0&0&\cdots&J_n
\end{pmatrix},$$ in which $J_i=\begin{pmatrix}
\lambda_i&1&0&\cdots&0\\
0&\lambda_i&1&\cdots&0\\
\vdots&\vdots&\ddots&\ddots&\vdots\\
0&0&\cdots&\lambda_i&1\\
0&0&\cdots&0&\lambda_i
\end{pmatrix}$ and $i=1,2,\cdots, n.$
Thus, $C-C^3$ is nilpotent $\Leftrightarrow P^{-1}CP-(P^{-1}CP)^3$ is nilpotent $\Leftrightarrow J_i-(J_i)^3$ is nilpotent.
We deduce $$J_i-(J_i)^3=\begin{pmatrix}
\lambda_i-(\lambda_i)^3&*&*&*\\
0&\lambda_i-(\lambda_i)^3&*&*\\
\vdots&\vdots&\ddots&\vdots\\
0&0&\lambda_{i-1}-(\lambda_{i-1})^3&*\\
0&0&0&\lambda_i-(\lambda_i)^3
\end{pmatrix}.$$
Then $J_i-(J_i)^3$ is nilpotent $\Leftrightarrow \lambda_i-(\lambda_i)^3=0 \Leftrightarrow \lambda_i(\lambda_i+1)(\lambda_i-1)=0 \Leftrightarrow \lambda_i=0$ or $1$ or $-1.$
\end{proof}

Now we come to present the explicit result for certain complex matrices.

\begin{thm}
Let $A,B\in {\Bbb C}^{n\times n}$. If $B^DA=0, BAB^{\pi}=0$ and the eigenvalues of $A$ and $B$ are $0$ or $1$, then
$$M=\begin{pmatrix}
A&I\\
B&0
\end{pmatrix}
\in ({\Bbb C}^{2n\times 2n})^{H}.$$ 
\end{thm}

\begin{proof}
As a result of Lemma 3.4, $M=\begin{pmatrix}
A&I\\
B&0
\end{pmatrix}
\in ({\Bbb C}^{2n\times 2n})^{H}$ can be easily obtained.
\end{proof}

\begin{exam}
Let $A=\begin{pmatrix}
1&1\\
0&0
\end{pmatrix}, B=\begin{pmatrix}
0&1\\
0&1
\end{pmatrix}\in \mathbb C^{2\times 2}.$ Then $B^DA=0, BAB^{\pi}=0$ and $A,B\in \mathcal A^{sD}.$  Then $$M=\begin{pmatrix}
A&I\\
B&0
\end{pmatrix}=\begin{pmatrix}
1&1&1&0\\
0&0&0&1\\
0&1&0&0\\
0&1&0&0
\end{pmatrix}.$$
$$\left|{\lambda}I-M\right|=\left|
\begin{matrix}
\lambda-1&-1&-1&0\\
0&\lambda&0&-1\\
0&-1&\lambda&0\\
0&-1&0&\lambda
\end{matrix}\right|=\lambda(\lambda+1)(\lambda-1)^2,$$ so $\lambda=0,1,-1.$ By Lemma 3.4, we have
$$
M=\begin{pmatrix}
A&I\\
B&0
\end{pmatrix}
\in ({\Bbb C}^{4\times 4})^{H}.$$
\end{exam}

\section{Communicative-like conditions}

In this section we consider the Hirano invertibility of anti-triangular matrices under certain communicative-like conditions.

\begin{thm}
Let $a\in {\mathcal A}^H, b\in \mathcal A^{sD}.$ If $b^Da=0, abb^{\pi}=bab^{\pi},$ then
$x=\begin{pmatrix}
a&1\\
b&0
\end{pmatrix}\in M_2({\mathcal A})^{H}.$
\end{thm}

\begin{proof}
Since $b$ is strongly Drazin invertible, it has Drazin inverse. Then $b$ can be written as $b=b_0+b_{00},$ where $b^s={b_0}^s\oplus 0,
b^{\pi}=0\oplus 1$ and $b_0$ is invertible. Since $b^Da=0$, $a$ has the form $a=\begin{pmatrix}
a_0&0\\
a_1&a_{00}
\end{pmatrix}$ and $a_{00}\in \mathcal A^{H}.$ From 
$$
\begin{pmatrix}
a&1\\
b&0
\end{pmatrix}=\begin{pmatrix}
a_0&0&1&0\\
a_1&a_{00}&0&1\\
b_0&0&0&0\\
0&b_{00}&0&0
\end{pmatrix}\hookrightarrow_S \begin{pmatrix}
a_0&1&0&0\\
b_0&0&0&0\\
a_1&0&a_{00}&1\\
0&0&b_{00}&0
\end{pmatrix},$$
where $S=E_{11}+E_{32}+E_{23}+E_{44}$ and $E_{ij}$ be the $4\times 4$ matrix with $(i,j)$  element equal to $1$ and others zero.

One easily checks that $a_0=bb^Dabb^D=bb^Da=0, b_0=b^2b^D\in \mathcal{A}^{sD}$. Moreover, we have $$b_0a_0b_0^{\pi}=b^2b^Dabb^D[1-b^2b^D(b^D)]=0.$$
In view of Theorem 3.1, we prove that $\begin{pmatrix}
a_0&1\\
b_0&0
\end{pmatrix}=\begin{pmatrix}
0&1\\
b_0&0
\end{pmatrix}$ has Hirano invertible.

Then we know that $\begin{pmatrix}
a&1\\
b&0
\end{pmatrix}$ is Hirano invertible if and only if $\begin{pmatrix}
a_{00}&1\\
b_{00}&0
\end{pmatrix}$ is Hirano invertible. From $abb^{\pi}=bab^{\pi}$ we get $a_{00}b_{00}=b_{00}a_{00}.$ So $a_{00}$ and $b_{00}$ can be written as $$a_{00}=a_3\oplus a_4, b_{00}=b_3\oplus b_4,$$ where $a_3$ is invertible, $a_4, b_3, b_4\in N(\mathcal A)$ and $a_ib_i=b_ia_i (i=3,4).$
So $$\begin{pmatrix}
a_{00}&1\\
b_{00}&0
\end{pmatrix}=\begin{pmatrix}
a_3&0&1&0\\
0&a_4&0&1\\
b_3&0&0&0\\
0&b_4&0&0
\end{pmatrix}\hookrightarrow_S \begin{pmatrix}
a_3&1&0&0\\
b_3&0&0&0\\
0&0&a_4&1\\
0&0&b_4&0
\end{pmatrix}.$$ 

Hence $\begin{pmatrix}
a_{00}&1\\
b_{00}&0
\end{pmatrix}$ is Hirano invertible if and only if $\begin{pmatrix}
a_i&1\\
b_i&0
\end{pmatrix}, i=3,4$ are Hirano invertible. Since $a_4,b_4\in \mathcal{A}$ are nilpotent and $a_4b_4=b_4a_4$, it follows by [2, Lemma 1] that $\begin{pmatrix}
a_4&1\\
b_4&0
\end{pmatrix}$ is nilpotent. 

Since $a_3$ is invertible, $b_3\in \mathcal{A}$ is nilpotent with $a_3b_3=b_3a_3,$
there exists $x_3\in N(\mathcal A)$ such that $a_3x_3+x_3^2=b_3$ and $x_3b_3\in \mathcal{A}, a_3x_3=x_3a_3, a_3+x_3$ is invertible . Then $\begin{pmatrix}
a_3&1\\
b_3&0
\end{pmatrix}\hookrightarrow_{S_0} \begin{pmatrix}
a_3+x_3&1\\
0&-x_3
\end{pmatrix},$ where $S_0=\begin{pmatrix}
1&0\\
-x_3&1
\end{pmatrix}.$ Following we claim that $\begin{pmatrix}
a_3+x_3&1\\
0&-x_3
\end{pmatrix}$ is Hirano invertible. 

Since $a_{00}=\begin{pmatrix}
a_3&0\\
0&a_4
\end{pmatrix}\in {\mathcal A}^H, a_4$ is nilpotent, therefore $a_3\in {\mathcal A}^H.$ According to Theorem 2.1, $a_3$ is the sum of an tripotent and a nilpotent. Addictively, $x_3$ is nilpotent, thus $a_3+x_3$ also can be written as the sum of an tripotent and a nilpotent. Therefore, $a_3+x_3\in {\mathcal A}^H.$

Then, $\begin{pmatrix}
a_3&1\\
b_3&0
\end{pmatrix}\in {\mathcal A}^H.$ We complete the proof.
\end{proof}

\begin{cor} Let $a\in {\mathcal A}^H, bc\in \mathcal{A}^{sD}.$ If $acb=cba, (cb)^Da=0,$ then
$\begin{pmatrix}
a&c\\
b&0
\end{pmatrix}\in M_2(\mathcal{A})^{H}.$
\end{cor}

\begin{proof}
We easily verify $\begin{pmatrix}
a&c\\
b&0
\end{pmatrix}=\begin{pmatrix}
1&0\\
0&c
\end{pmatrix}\begin{pmatrix}
a&1\\
b&0
\end{pmatrix}$ and $\begin{pmatrix}
a&1\\
b&0
\end{pmatrix}\begin{pmatrix}
1&0\\
0&c
\end{pmatrix}=\begin{pmatrix}
a&1\\
cb&0
\end{pmatrix}.$ In light of Lemma 2.3, $\begin{pmatrix}
a&c\\
b&0
\end{pmatrix}$ is Hirano invertible if and only if $\begin{pmatrix}
a&1\\
cb&0
\end{pmatrix}$ is Hirano invertible. Since $a, cb\in (\mathcal A)^{sD}, acb=cba,$ apply Theorem 4.1 to $\begin{pmatrix}
a&1\\
cb&0
\end{pmatrix},$ we have $\begin{pmatrix}
a&1\\
cb&0
\end{pmatrix}\in M_2(\mathcal A)^{H}.$
Accordingly, $\begin{pmatrix}
a&b\\
c&0
\end{pmatrix}\in M_2(\mathcal A)^{H}.$
\end{proof}

\begin{cor}
Let $a,b\in \mathcal A^{sD}.$ If $b^Da=0, ab=ba,$ then $\begin{pmatrix}
a&1\\
b&0
\end{pmatrix}\in M_2({\mathcal A})^{H}.$
\end{cor}

\begin{proof}
This is obvious by Corollary 4.2.
\end{proof}

\begin{thm}
Let $a\in {\mathcal A}^H, b\in \mathcal A^{sD}.$ If $b=ba,$ then
$\begin{pmatrix}
a&1\\
b&0
\end{pmatrix}\in M_2({\mathcal A})^{H}.$
\end{thm}

\begin{proof}
When $b$ can be written as $\begin{pmatrix}
b_1&0\\
0&b_2
\end{pmatrix},$ $b_1$ is invertible, ${b_2}^t=0$ for some $t\in \mathbb N.$ Since $b=ba,$ we have $a=\begin{pmatrix}
bb^D&0\\
a_3&a_2
\end{pmatrix},$ and $b_2=b_2a_2.$ Then $$\begin{pmatrix}
a&1\\
b&0
\end{pmatrix}=\begin{pmatrix}
bb^D&0&1&0\\
a_3&a_2&0&1\\
b_1&0&0&0\\
0&b_2&0&0
\end{pmatrix}\hookrightarrow_P \begin{pmatrix}
bb^D&1&0&0\\
b_1&0&0&0\\
a_3&0&a_2&1\\
0&0&b_2&0
\end{pmatrix}.$$
Notice that from $b_2=b_2a_2,$ we can obtain $b_2=\begin{pmatrix}
b_{21}&0\\
b_{23}&0
\end{pmatrix},$ and $b_{21}=b_{21}a_{21}.$ Then $$\begin{pmatrix}
a_2&1\\
b_2&0
\end{pmatrix}=\begin{pmatrix}
a_{21}&0&1&0\\
0&a_{22}&0&1\\
b_{21}&0&0&0\\
b_{23}&0&0&0
\end{pmatrix}\hookrightarrow \begin{pmatrix}
a_{21}&1&0&0\\
b_{21}&0&0&0\\
0&0&a_{22}&1\\
b_{23}&0&0&0
\end{pmatrix}.$$
Because ${a_{22}}^n=0$ for some $n\in \mathbb N,$ so does $\begin{pmatrix}
a_{22}&1\\
0&0
\end{pmatrix}.$ Therefore, $\begin{pmatrix}
a_2&1\\
b_2&0
\end{pmatrix}$ Hirano invertible is equal to $\begin{pmatrix}
a_{21}&1\\
b_{21}&0
\end{pmatrix}$ Hirano invertible. Since ${b_2}^t=0\Rightarrow {b_{21}}^t=0, b_{21}=b_{21}a_{21}, a_{21}$ is invertible, there exists a nilpotent operator $x_{21}$ such that $x_{21}^2+a_{21}x_{21}=b_{21},$ and $a_{21}+x_{21}$ is invertible, Thus, $$\begin{pmatrix}
a_{21}&1\\
b_{21}&0
\end{pmatrix}\sim \begin{pmatrix}
a_{21}+x_{21}&1\\
0&-x_{21}
\end{pmatrix}.$$
According to Theorem 4.1 ,$\begin{pmatrix}
a_{21}+x_{21}&1\\
0&-x_{21}
\end{pmatrix}$ is Hirano invertible, so does $\begin{pmatrix}
a_{21}&1\\
b_{21}&0
\end{pmatrix}.$ Thus, $\begin{pmatrix}
a&1\\
b&0
\end{pmatrix}\in M_2({\mathcal A})^{H}.$
\end{proof}

\begin{cor}
Let $a\in {\mathcal A}^H, b\in \mathcal A^{sD}.$ If $b=ab,$ then
$\begin{pmatrix}
a&1\\
b&0
\end{pmatrix}\in M_2({\mathcal A})^{H}.$
\end{cor}

\begin{proof}
When $b$ can be written as $\begin{pmatrix}
b_1&0\\
0&b_2
\end{pmatrix},$ $b_1$ is invertible, ${b_2}^t=0$ for some $t\in \mathbb N.$ Since $b=ab,$ we have $a=\begin{pmatrix}
bb^D&a_3\\
0&a_2
\end{pmatrix},$ and $b_2=a_2b_2.$ Then $$\begin{pmatrix}
a&1\\
b&0
\end{pmatrix}=\begin{pmatrix}
bb^D&a_3&1&0\\
0&a_2&0&1\\
b_1&0&0&0\\
0&b_2&0&0
\end{pmatrix}\hookrightarrow_P \begin{pmatrix}
bb^D&1&a_3&0\\
b_1&0&0&0\\
0&0&a_2&1\\
0&0&b_2&0
\end{pmatrix}.$$
Notice that from $b_2=a_2b_2,$ we can obtain $b_2=\begin{pmatrix}
b_{21}&b_{22}\\
0&0
\end{pmatrix},$ and $b_{21}=a_{21}b_{21}.$ Then $$\begin{pmatrix}
a_2&1\\
b_2&0
\end{pmatrix}=\begin{pmatrix}
a_{21}&0&1&0\\
0&a_{22}&0&1\\
b_{21}&b_{22}&0&0\\
0&0&0&0
\end{pmatrix}\hookrightarrow \begin{pmatrix}
a_{21}&1&0&0\\
b_{21}&0&b_{22}&0\\
0&0&a_{22}&1\\
0&0&0&0
\end{pmatrix}.$$
Because ${a_{22}}^n=0$ for some $n\in \mathbb N,$ so does $\begin{pmatrix}
a_{22}&1\\
0&0
\end{pmatrix}.$ Therefore, $\begin{pmatrix}
a_2&1\\
b_2&0
\end{pmatrix}$ Hirano invertible is equal to $\begin{pmatrix}
a_{21}&1\\
b_{21}&0
\end{pmatrix}$ Hirano invertible. Since ${b_2}^t=0\Rightarrow {b_{21}}^t=0, b_{21}=a_{21}b_{21}, a_{21}$ is invertible, there exists a nilpotent operator $x_{21}$ such that $x_{21}^2+a_{21}x_{21}=b_{21},$ and $a_{21}+x_{21}$ is invertible, Thus, $$\begin{pmatrix}
a_{21}&1\\
b_{21}&0
\end{pmatrix}\sim \begin{pmatrix}
a_{21}+x_{21}&1\\
0&-x_{21}
\end{pmatrix}.$$
According to Theorem 4.4 ,$\begin{pmatrix}
a_{21}+x_{21}&1\\
0&-x_{21}
\end{pmatrix}$ is Hirano invertible, so does $\begin{pmatrix}
a_{21}&1\\
b_{21}&0
\end{pmatrix}.$ Thus, $\begin{pmatrix}
a&1\\
b&0
\end{pmatrix}\in M_2({\mathcal A})^{H}.$
\end{proof}

\section{Numeral examples}

In this section, three specific examples would be given to illustrate the main results.

\begin{exam}
$A=\begin{pmatrix}
0&0&1\\
1&0&1\\
1&0&0
\end{pmatrix},$ we can judge $A\in M_3(A)^H.$
\end{exam}

\begin{proof}
$A=\begin{pmatrix}
0&0&1\\
1&0&1\\
1&0&0
\end{pmatrix}=\begin{pmatrix}
0&0&1\\
1&0&0\\
1&0&0
\end{pmatrix}+\begin{pmatrix}
0&0&0\\
0&0&1\\
0&0&0
\end{pmatrix}=B+C, B^3=B, C\in N(\mathcal A).$ According to Theorem 2.1, $\begin{pmatrix}
0&0&1\\
1&0&1\\
1&0&0
\end{pmatrix}\in M_3(A)^H.$
\end{proof}

\begin{exam}
Let $A=\begin{pmatrix}
0&0\\
0&1
\end{pmatrix}\in {\mathcal A}^H, B=\begin{pmatrix}
1&0\\
1&0
\end{pmatrix}\in {\mathcal A}^{sD}.$ Then $x=\begin{pmatrix}
A&I\\
B&0
\end{pmatrix}\in M_4(\mathcal A)^H.$
\end{exam}

\begin{proof}
We have $B^DA=0, BAB^{\pi}=0.$ Using Theorem 3.1, $x=\begin{pmatrix}
A&I\\
B&0
\end{pmatrix}=\begin{pmatrix}
0&0&1&0\\
0&1&0&1\\
1&0&0&0\\
1&0&0&0
\end{pmatrix}\in M_4(\mathcal A)^H$ obtained.
\end{proof}

\begin{exam}
Let $A=\begin{pmatrix}
0&0\\
0&1
\end{pmatrix}\in {\mathcal A}^H, B=\begin{pmatrix}
1&0\\
-1&0
\end{pmatrix}\in {\mathcal A}^{sD}.$ Then $x=\begin{pmatrix}
A&I\\
B&0
\end{pmatrix}\in M_4(\mathcal A)^H.$
\end{exam}

\begin{proof}
We have $B^DA=0, BAB^{\pi}=ABB^{\pi}.$ Using Theorem 4.1, $x=\begin{pmatrix}
A&I\\
B&0
\end{pmatrix}=\begin{pmatrix}
0&0&1&0\\
0&1&0&1\\
1&0&0&0\\
-1&0&0&0
\end{pmatrix}\in M_4(\mathcal A)^H$ obtained.
\end{proof}

\end{document}